\documentclass[11pt]{amsart}

\newcommand{\om}{\omega}
\newcommand{\IM}{\mathrm{Im}}
\newcommand{\Sl}{\mathrm{sl}}

\begin{document}

\title[Discrete series of representations
for the modular double]
{Discrete series of representations
for the modular double of
 $U_q(\Sl(2,\mathbb{R}))$}

\author{L. D. Faddeev}

\address{Saint Petersburg Department of
Steklov Mathematical Institute}

\begin{abstract}
Modular double of quantum group $U_q (\Sl(2))$ with deformation parameter
$q=e^{i\pi\tau}$ is a natural object explicitly taking into account
the duality $\tau \mapsto 1/\tau$. The use of the modular double in
CFT allows to consider the region $1<c<25$ for the central charge of the
Virasoro algebra when $|\tau|=1$. In this paper a new discrete series
of representations for the modular double of $U_q (\Sl(2,\mathbb{R}))$
is found for such $\tau$.
\end{abstract}

\maketitle

\section{Introduction}

Modular double of quantum group $U_q(\Sl(2))$ is loosely speaking a
tensor product of two quantum groups with deformation parameters
$q=e^{i\pi\tau}$ and $\tilde q=e^{i\pi/\tau}$. It was introduced by
me in \cite{1} as an outcome of my quest for regularization of universal
$R$-matrix at roots of unity. Its relevance to CFT, especially ones
with noncompact field variables such as Liouville model, is rather
evident. In several papers \cite{2}, \cite{3}, \cite{4},
\cite{5a} for the noncompact  real form $U_q(\Sl(2,\mathbb{R}))$
 the continuous series of unitary representations, parameterized
by a positive real number (spin), was introduced.
An interesting feature is that these representations do not decompose into
tensor product of representations for dual factors. These representations
exist for two regimes of deformation parameter $\tau$.

Regime I: $\tau$ real, positive.

Regime II: $|\tau|=1$, $\IM \tau\ge 0$.

In this note we present one more series of representations for
regime II which is discrete: the representations will be
parameterized by two integers: a positive $n$ and an integer $m$
with values $m=1, \ldots,n$.

In the first section we remind known results, mentioned above, and
in the second section we present new representations.

In the course of exposition the role of a special function, which I
call noncompact dilogarithm, will become clear. It is my conviction
that this function, properties of which are lucidly explained in
\cite{5}, is the main technical object of CFT and integrable  models, as
well as quantum Teichmuller theory, for recent example see \cite{6}, \cite{7}.

This work is party supported by RFBR grant 08-01-00638 and programme
``Mathematical problems of nonlinear dynamics'' of Russian Academy
of Sciences.  I am grateful to P.~Kulish and M.~Semenov-Tian-Shansky
for valuable discussions.

\section{Warm up---continuous series}

The set of generators for modular double of $U_q(\Sl(2))$ consists
of two triples $E,F,K$ and $\tilde E,\tilde F,\tilde K$ with
relations
\begin{align}
KE&=q^2EK
\\
KF&=q^{-2}FK
\\
EF-FE&=\frac{1}{q-q^{-1}}(K-K^{-1})
\end{align}
with complex deformation parameter $q$,
$$
q=e^{i\pi \tau}
$$
and similar relation for $\tilde E,\tilde F,\tilde K$ with parameter $\tilde q$,
$$
\tilde q=e^{i\pi/\tau}
$$
Note that  use of $1/\tau$ instead of traditional $-1/\tau$ has several conveniences, but is not principal.

Generators $E,F,K$ and $\tilde E,\tilde F,\tilde K$ mutually commute.

We shall use expression for these generators via simpler variables,
in analogue with  Gelfand--Kirillov variables for the non-deformed
universal enveloping algebra $U(\Sl(2))$. Namely consider the Weyl pair
$u,v$ with relation
$$
uv=q^2vu
$$
and central element $Z$. Then
$$
E=i\frac{v+u^{-1}Z}{q-q^{-1}},
\qquad
F=i\frac{u+v^{-1}Z^{-1}}{q-q^{-1}},
\qquad
K=q^{-1}uv
$$
satisfy relations (1)--(3). The traditional Casimir
$$
C=qK+q^{-1}K^{-1}+(q-q^{-1})^2FE
$$
is expressed via $Z$ as
$$
C=-(Z+Z^{-1}).
$$
In \cite{3}, \cite{4} and \cite{5a} similar and equivalent parameterizations
were used.

The dual generators are expressed similarly via
$\tilde u,\tilde v,\tilde Z$, where
$$
\tilde u=u^{1/\tau}, \qquad
\tilde v=v^{1/\tau}, \qquad
\tilde Z=Z^{1/\tau}.
$$
The commutativity between, $E,F,K$ and $\tilde E,\tilde F,\tilde K$
is due to Euler formula
$$
e^{2\pi i}=1.
$$
The representations, introduced in \cite{3}, \cite{4} are realized in
$L_2(\mathbb{R})$ with the Weyl pairs expressed as usual via
multiplication and shift. To be more explicit, introduce convenient
Weierstrass-like complex parameters $\om, \om'$ with a constraint
$$
\om\om'=-\frac{1}{4},
$$
through which $\tau$ is expressed as
$$
\tau=\frac{\om'}{\om},
$$
so that turning to dual means interchange of $\om$ and $\om'$.

The elements of $L_2(\mathbb{R})$ are functions $f(x)$ of real
variable $x$ on which $u,v$ act as follows
\begin{align}
uf(x)&=  e^{-i\pi x/\om}f(x),
\\
vf(x)&= f(x+2\om').
\end{align}
The shift $x \mapsto x+2\om'$ leaves the real axis, and function
$e^{-i\pi x/\om}$ is unbounded, so the domain for operators $u$ and
$v$ must be chosen accordingly. A  suitable variant is to use dense
domain, consisting of functions $f(x)$, analytic in the whole plane
$\mathbb{C}$ and rapidly vanishing along the contour, parallel to
real axis. Functions like
$$
f(x)=e^{-x^2}P(x),
$$
where $P(x)$ is a polynomial, are typical examples.

With this convention it is easy to calculate the operator, adjoint
to $v$ with respect to scalar product of $L_2(\mathbb{R})$
\begin{equation}
(f,g)=\int^{\infty}_{-\infty}\overline{f(x)}g(x)dx.
\end{equation}
Indeed
\begin{align*}
(f,vg)&
=\int^{\infty}_{-\infty}\overline{f(x)}g(x+2\om')dx
\\
&=\int^{\infty}_{-\infty}\bar f(x)g(x+2\om')dx
\\
&=\int_{-\infty+2\om'}^{\infty+2\om'}\bar
f(x-2\om')g(x)dx
\\
&=\int^{\infty}_{-\infty}
\overline{f(x-\overline{2\om'})}g(x)dx
\\
&=(v^*g,f),
\end{align*}
where
$$
v^*f(x)=f(x-\overline{2\om'}).
$$
Here we denoted by $\bar f(x)$ an analytic function of $x$ with
coefficients in expansion in powers of $x$ being complex conjugate
to those of $f$.

These considerations allow us to introduce the real form for our
algebra. In fact there are two possibilities, which we shall call
regime I and regime II.

Regime I: $\tau$ real positive, $\om,\om'$ pure imaginary
$$
\bar \om=-\om,\qquad \bar \om'=-\om'.
$$
Operators $u,v$ and $\tilde u$, $\tilde v$, are Hermitian
and positive definite $u^*=u$, $v^*=v$,
$\tilde u^*=\tilde u$,
$\tilde v^*=\tilde v$ and the same
is true for $E,F,K$ and $\tilde E,\tilde F,\tilde K$, if $Z$
 is Hermitian. In irreducible representation $Z$ must be  a real number.
We shall parameterize it as
\begin{equation}
Z=e^{i\pi a/\om}, \qquad \tilde Z=e^{i\pi a/\om'}
\end{equation}
with real $a$---and naturally call $a$ the spin.

More mathematical details we can find in \cite{3}, \cite{4}, \cite{5a}.

Regime II. $|\tau|=1$
$$
\bar \om=-\om'.
$$
Operators $u,v$ are hermit conjugate to $\tilde u,\tilde v$
$$
u^*=\tilde u,\qquad v^*=\tilde v
$$
and the same is true for numbers $Z$ and $\tilde Z$
$$
\overline
{\tilde{Z}}=Z,
$$
if we use the same parametrization (7) with real $a$.

Let us note that the third Weierstrass parameter
$$
\om''=\om+\om'
$$
is pure imaginary in both cases.

This concludes our reminder. It clear that the scalar product (6)
was exhaustively used in this presentation. In the next section we
shall explore another possibility.

\section{The discrete series}

We shall introduce here another scalar product for function $f$, on
which the Weyl operators act.  Now they are analytic functions
$f(z)$ and action of $u,v,\tilde u,\tilde v$ are the same as in (4),
(5) with change $x\mapsto z$. The scalar product will be
\begin{equation}
\langle f,g\rangle =\int_D\overline{f(z)}
S(\bar z,z)g(z)\frac{dzd\bar z}{2\pi i}
\end{equation}
with integration over domain $D$ in $\mathbb{C}$. The positive
kernel $S(\bar z,z)$ should be found from the condition of
hermiticity. We shall find, that reasonable solution exists for
discrete values of spin
$ a $, namely
$ a = -n\om'' $
with positive integer
$ n $
in regime II, so that spin
$ a $
is pure imaginary and the following relation holds:
$$
\bar Z=\tilde Z^{-1}.
$$
Begin with relation
$$
\tilde K^*=K.
$$
We have
\begin{align*}
\langle \tilde Kf,g\rangle
&=\int \overline{\tilde q^{-1}e^{-i\pi z/\om'}f(z+2\om)}
S(\bar z,z)g(z)\frac{dzd\bar z}{2\pi i}
\\
&=q \int \bar f(\bar
z-2\om')S(\bar z,z)e^{-i\pi \bar z/\om}g(z)\frac{dzd\bar z}{2\pi i}
\\
&=q^{-1}\int  \overline{f(z)}S(\bar z+2\om',z)e^{-i\pi \bar z/\om}
g(z)\frac{dzd\bar z}{2\pi i}
\end{align*}
and
\begin{align*}
\langle f,Kg\rangle
&= \int  \overline{f(z)}S(\bar z,z)q^{-1}e^{-i\pi
z/\om}g(z+2\om')\frac{dzd\bar z}{2\pi i}
\\
&=q\int \overline{f(z)} S(\bar z,z-2\om')e^{-i\pi z/\om}g(z)\frac{dzd\bar z}{2\pi i}.
\end{align*}
The equation
$$
qS(\bar z,z-2\om')e^{-i\pi z/\om}=q^{-1}S(\bar z+2\om',z)e^{-i\pi \bar z/\om}
$$
is satisfied by the kernel $S(\bar z,z)$ of the form
\begin{equation}
S(\bar z,z)=\exp i\pi \big(z^2-\bar z^2\big)\Phi(\bar z-z).
\end{equation}
Now turn to next condition $E^*=\tilde E$. We have
\begin{align*}
\langle f,Eg\rangle
&=\frac{i}{q-q^{-1}}\int\overline{f(z)}S(\bar z,z)
\left(g(z+2\om')+Ze^{i\pi z/\om}g(z)\right)\frac{dzd\bar z}{2\pi i}
\\
&=\frac{i}{q-q^{-1}}\int\overline{f(z)}
\left(S(\bar z,z-2\om')
+ZS(\bar z,z)e^{i\pi z/\om}\right)g(z)\frac{d zd\bar z}{2\pi i}
\end{align*}
and
\begin{align*}
(\tilde Ef,q)
&=\int\overline{\frac{i}{\tilde q-\tilde q^{-1}}
\left(f(z+2\om)
+Z e^{i\pi z/\om'}f(z)\right)}
S(\bar z,z)g(z)\frac{d zd\bar z}{2\pi i}
\\
&=\frac{i}{q-q^{-1}}\int\overline{f(z)}\left(S(\bar z+2\om',z)
+\overline{\tilde Z} S(\bar z,z)e^{i\pi \bar z/\om}\right)g(z)
\frac{dzd\bar z}{2\pi i},
\end{align*}
leading to the equation
$$
S(\bar z,z-2\om')+Ze^{i\pi z/\om}S(\bar z,z)=S(\bar z+2\om',z)+Z^{-1}e^{i\pi
\bar z/\om}S(\bar z,z).
$$
Introducing (9) here leads to the following equation for $\Phi(t)$, where
$t=\bar z-z$,
$$
\left(q^{-1}e^{-i\pi t/\om}-q\right)\Phi(t+2\om')
=\left(Z^{-1}-Ze^{-i\pi t/\om}\right)\Phi(t).
$$

    Substituting here
$$
    Z = e^{-i\pi a/\om}
$$
    (change of sign will be convenient in what follows) and noting, that
$$
    q = e^{i\pi \frac{\om'}{\om}} = - e^{i\pi\frac{\om''}{\om}}
$$
    we rewrite the last equation as follows
\begin{equation}
\label{Peq1}
    \frac{\Phi(t+2\om')}{\Phi(t)} = \frac{\sin \pi \frac{t+2a}{2\om}}{
    \sin \pi \frac{t+2\om''}{2\om}} .
\end{equation}

    The condition
$$
    (\tilde{F}f, q) = (f, Fq)
$$
    leads to the same equation. The conditions
$$
    (Ef,q) = (f,\tilde{E}q) , \quad (Ff,q) = (f,\tilde{F}q)
$$
    give equation with interchange of
$ \om $ and $ \om' $
\begin{equation}
\label{Peq2}
    \frac{\Phi(t+2\om)}{\Phi(t)} = \frac{\sin \pi \frac{t+2a}{2\om'}}{
    \sin \pi \frac{t+2\om''}{2\om'}} .
\end{equation}
    Both equations
(\ref{Peq1}) and (\ref{Peq2})
    fix the solution
$ \Phi(t) $
    essentially uniquelly up to a constant factor.

    We should look for solution
$ \Phi(t) $,
    which is positive for pure imaginary
$ t $.
    It is easy to see, that for the discrete values of ``spin''
$ a $
\begin{equation}
\label{spin}
    a = n\om'' , \quad n=1,2,\ldots
\end{equation}
    such solution is given by
$ \Phi(t)=1 $ for $ n=1 $
    and
\begin{equation}
\label{Phit}
    \Phi(t) = \prod_{m=1}^{n-1} \sin \pi \Big(\frac{t+2m\om''}{2\om'}\Big)
    \sin \pi \Big(\frac{t+2m\om''}{2\om}\Big) ,
\end{equation}
    which can be rewritten as
$$
    \Phi(t) = \prod_{m=1}^{n-1} \Big| \sin \pi \Big(\frac{t+2m\om''}{2\om'}\Big)
    \Big|^{2}
$$
    for imaginary
$ t $,
    so that it is manifestly nonnegative.

The lines of zeros of $\Phi$
$$
t=-2m\om''
$$
divide all complex plane into $n$ parts, including two half-planes
$y<\mu$, $y>\mu n$
and $n-2$ strips of finite width
$$
    \mu m < y < \mu (m+1) , \quad m=1, \ldots n-1 ,
$$
where
$$
y=\frac{z-\overline z}{2i}, \qquad \mu=\frac {\omega''}{i},
$$
giving us $n$ domains of integration in (8). I conjecture that the
corresponding representations are nonequivalent. If it is true, we
acquire the series of representations, parameterized by discrete
spin $a$ from
(\ref{spin})
    and additional number $m=1,\dots ,n-1$.

    For generic value of imaginary spin
$ a $
    the solution of equations
(\ref{Peq1}) and
(\ref{Peq2})
    can be expressed via noncompact dilogarithm
$ \gamma(t) $,
    mentioned in the introduction.
    Function
$ \gamma(z) $
    satisfies dual equations
\begin{align}
\label{d1}
    \frac{\gamma(\zeta+\om')}{\gamma(\zeta-\om')} =& 1 + e^{-i\pi\zeta/\om} ,\\
\label{d2}
    \frac{\gamma(\zeta+\om)}{\gamma(\zeta-\om)} =& 1 + e^{-i\pi\zeta/\om'}
\end{align}
    and solution
$ \Phi(t) $
    can be written as
\begin{equation}
\label{Phim}
    \Phi(t) = e^{2\pi i(a+\om'')} \frac{\gamma(t-\om''+2a)}{\gamma(t+\om'')} .
\end{equation}
    Of course for
$ a=n\om'' $
    it reduces to
(\ref{Phit})
    by means of relation
$$
    \frac{\gamma(\zeta+\om'')}{\gamma(\zeta-\om'')} = - 4 e^{2\pi i\zeta\om''}
    \sin \frac{\pi\zeta}{2\om'} \sin \frac{\pi\zeta}{2\om} ,
$$
    which follows from
(\ref{d1}), (\ref{d2}).
    I believe, that
(\ref{Phim})
    reduces to function positive for imaginary
$ t $
    only in the case considered above.

I cannot help mentioning the correspondence between parameters
$ n,m $
and zeros of function $\gamma(z)$ on the upper half-plane: These
zeros are neatly described by $z=\om''+2p\om+2 q\om'$, with $p,q$  being
nonnegative integers. It is clear that on each level $\IM z=n\om''$
there are $n$ zeros---in correspondence with the number of
representations.

This makes us to think about the deformed Plancherel formula.
However we shall leave it for the future work.

\section{Conclusion}

The representations, found here, are meaningful only in Regime II.
Indeed, if we repeat calculation for Regime I, we shall get the same
expression for the function $\Phi$, however the arguments in sines
will be real and function $\Phi$ will acquire infinite number of
zeros, dense in the case of irrational  $\tau$ and it will be
difficult to make any sense for this. The case of root of unity may
deserve further investigation.

I believe that Regime II is more interesting than Regime I both
mathematically and in the application to physics. Indeed, $\om,\om'$
define lattice in $\mathbb{C}$ in Regime II, whereas they  are
confined to the imaginary axis in Regime I. In the Liouville model
the central charge of Virasoro algebra is given by
$$
C=1+6\Big(\tau+\frac{1}{\tau}+2\Big)
$$
(see \cite{8}, \cite{9})
so that $c>25$ in Regime I and $1<c<25$ in Regime II.

Regime II,  is interesting for the application of Liouville model to
Polyakov noncritical string \cite{10} and new representations open way
to search for the quantization of Liouville model, different from
mere continuation from the case $c>25$ \cite{11} along the
continuous series of section 2.

This leaves us with the problem of adequate use of the representations,
found in this paper.


\begin{thebibliography}{**}

\bibitem{1}  L. D. Faddeev,
Modular double of a quantum group,
\textit{Math. Phys. Studies} \textbf{21} (2000), 149--156.

\bibitem{2} A. G. Bytsko and J. Teschner,
Quantization of models non-compact quantum group symmetry:
modular XXZ magnet and lattice sinh-Gordon model.
\textit{J. Phys. A: Math. Gen.} \textbf{39} (2006), 12927--12981.

\bibitem{3} B. Ponsot and J. Teschner,
Clebsch--Gordan and Racah--Wigner coefficients for a
continuous series of representations of $U_q(\mathrm{sl}(2,\mathbb{R}$)).
\textit{Commun. Math. Phys.} \textbf{224} (2001), 613--655.

\bibitem{4} S. Kharchev, D. Lebedev, and M. Semenov-Tian-Shansky,
Unitary  representations of $U_q(\textrm{sl}(2,\mathbb{R}$)),
the modular double, and the multiparticle $q$-deformed Toda chains.
\textit{Commun. Math. Phys.} \textbf{225} (2002), 573--609.

\bibitem{5a} K.~Schmudgen,
Operator representation of $ U_{q}(\textrm{sl}(2))$,
\textit{Lett. Math. Phys.} \textbf{37} (1996), 211--222.

\bibitem{5} A. Yu. Volkov,
Noncommutative hypergeometry. \textit{Commun. Math. Phys.}
\textbf{258} (2005), 257--273.

\bibitem{6} R. M. Kashaev,
Quantization of Teichm\"uller spaces and the quantum dilogarithm.
\textit{Lett. Math. Phys.} \textbf{43}, no.~2 (1998), 105--115.

\bibitem{7} V. V. Fock and L. O. Chekhov,
Quantum Teichm\"uller spaces. [In Russian]
\textit{Teoret. Mat. Fiz.}, \textbf{120}, no.~3 (1999), 511--528;
translation in \textit{Theor. Math. Phys.} \textbf{120},
(1999), 1245--1259.

\bibitem{8} J.-L. Gervais and  A. Neveu,
The dual string spectrum in Polyakov's quantization (I).
\textit{Nucl. Phus. B} \textbf{199} (1982), 59--76.

\bibitem{9} T. L. Curtright and C. B. Thorn,
Conformally invariant quantization of the Liouville theory.
\textit{Phys. Rev. Lett.} \textbf{48} (1982), 1309--1313.

\bibitem{10} A. M. Polyakov,
Quantum geometry of bosonic strings.
\textit{Phys. Lett. B} \textbf{103} (1981), 207--210.

\bibitem{11} V. V. Bazhanov, S. L. Lukyanov, and A. B. Zamolodchikov,
Integrable structure of conformal field theory, quantum KdV
theory and thermodynamic Bethe ansatz.
\textit{Commun. Math. Phys.} \textbf{177} (1996), 381;
{\tt arxiv:hep-th/9412229}.

\end{thebibliography}
\end{document}